\documentclass[a4paper,12pt]{article}
\usepackage{amsmath,amsfonts,amssymb,amscd}
\usepackage{hyperref}
\input{xypic}
\usepackage{tikz}
\usepackage{multicol}

\def\s{\scriptstyle }
\def\a{\alpha}

\def\l{\lambda}

\def\d{\partial}

\def\x{\mathbf x}
\def\y{\mathbf y}

\def\mfS{{\mathfrak S}}
\def\Sym{\mathfrak{ S\hspace{-0.07 em}y\hspace{-0.07 em}m}}

\def\N{{\mathbb N}}
\def\Z{{\mathbb Z}}

\def\C{{\mathbb C}}

\def\cH{{\mathcal H}}

\def\cR{{\mathcal R}}

\def\cD{{\mathcal D}}

\def\tiret{\ \raise3.5pt\hbox{\vrule height 0.5pt depth 0pt width 20pt}\ }

\def\moins{\raise 1pt\hbox{{$\scriptstyle -$}}}
\def\plus{\raise 1pt\hbox{{$\scriptstyle +$}} }
\def\demi{{ \s \frac{ 1}{2}}}
\def\hfl#1{{\stackrel{#1}{\hbox to 30pt{\rightarrowfill}}}} 

\def\QED{\hfill QED}

\newtheorem{theorem}{Theorem}
\newtheorem{proposition}[theorem]{Proposition}
\newtheorem{lemma}[theorem]{Lemma}
\newtheorem{corollary}[theorem]{Corollary}

\def\proof{\noindent{\it Proof.}\ }

\begingroup
\catcode`\@=11
\gdef\petitematrice#1{\null\vcenter {\normalbaselines \m@th
\ialign {\hfil $##$\hfil &&\thinspace  \hfil $##$\hfil\crcr
\mathstrut \crcr \noalign {\kern -\baselineskip } #1\crcr
\mathstrut \crcr \noalign {\kern -\baselineskip }}}}

\gdef\moyennematrice#1{\null\vcenter {\normalbaselines \m@th
\ialign {\hfil $##$\hfil &&\   \hfil $##$\hfil\crcr
\mathstrut \crcr \noalign {\kern -\baselineskip } #1\crcr
\mathstrut \crcr \noalign {\kern -\baselineskip }}}}
\endgroup
\def\mymatrixserre#1{\left\lbrack\vcenter{\offinterlineskip
     \halign{\vrule height 6pt depth 2pt width 0pt
      \hfil$\scriptstyle##$\hfil&&\kern
       4pt\hfil$\scriptstyle##$\hfil \crcr#1\crcr}}\right\rbrack}

\newdimen\carresize
\newdimen\thickness
\def\CARRE#1{\hbox{\vrule width \thickness
   \vbox to \carresize{\hrule height \thickness\vss
      \hbox to \carresize{\hss#1\hss}
   \vss\hrule height\thickness}
\unskip\vrule width \thickness}
\kern-\thickness}
\def\vsquare#1{\vbox{\CARRE{$#1$}}\kern-\thickness}

\def\smallyoung#1{%
  \carresize=12pt%
  \thickness=0.5pt%
  \vcenter{%
    \vbox{\smallskip\offinterlineskip%
      \halign{&\vsquare{##}\cr #1}}}}

\def\young#1{%
  \carresize=16pt%
  \thickness=0.5pt%
  \vcenter{%
    \vbox{\smallskip\offinterlineskip%
      \halign{&\vsquare{##}\cr #1}}}}

\def\bigyoung#1{%
  \carresize=25pt%
  \thickness=0.7pt%
  \vcenter{%
    \vbox{\smallskip\offinterlineskip%
      \halign{&\vsquare{##}\cr #1}}}}

\newdimen\vcadre\vcadre=0.2cm 
\newdimen\hcadre\hcadre=0.2cm 
\def\GrTeXBox#1{\vbox{\vskip\vcadre\hbox{\hskip\hcadre%
      $#1$%
   \hskip\hcadre}\vskip\vcadre}}
\def\arx#1[#2]{\ifcase#1 \relax \or%
  \ar @{-}[#2]  \or%
  \ar @2{-}[#2] \or%
  \ar @{--}[#2] \or%
  \ar @2{.}[#2] \or%
  \ar @{~}[#2]  \fi}
\tikzstyle{Pas} = [color=blue!70,line width=2pt]
\tikzstyle{Pas2} = [color=red!90,line width=2pt]

\newdimen\unit
\def\o{$\scriptscriptstyle{{\rm o}}$}
\def\grape(#1,#2)#3{\raise#2\unit\rlap{\kern#1\unit #3}\ignorespaces}

\def\gg{{\unit=1mm
\hbox {\grape(2,1.2){'}
       \grape(1,2)\o
       \grape(2,2)\o
       \grape(3,2)\o
       \grape(1.5,1)\o
       \grape(2.5,1)\o
       \grape(2,0)\o\
}\kern 3.5 \unit}}

\def\gfill{\leaders\hbox to 1.2em{\hss\gg\hss}\hfill}
\def\frise{\centerline{\hbox to 8cm{\gfill} }\bigskip}

\def\Pf{{\mathfrak{Pf}}}
\def\M{{\mathfrak{M}}}
\def\a{\mathbf a}
\def\h{\mathbf h}
\def\q{\mathbf q}

\def\Coin{\mathfrak{M}^e}
\def\Bez{\mathcal B\hspace{-0.06 em} {\sl e} \hspace{-0.06 em}{\sl z}}

\def\KL{K\hspace{-0.35 em}L}
\def\LK{L\hspace{-0.25 em}K}
\date{}

\begin{document}

\title{\bf Hankel Pfaffians, Discriminants and 
  Kazhdan-Lusztig bases}

\author{Alain Lascoux}

\maketitle

\begin{abstract}
We use Kazhdan-Lusztig bases of representations of the 
symmetric group to express Pfaffians with entries 
$(a_i-a_j) h_{i+j}$. In the case where the parameters $a_i$
are specialized to successive powers of $q$, and the $h_i$
are complete functions, we obtain the $q$-discriminant. 
\end{abstract}

    \frise
     
Hankel matrices are matrices constant along anti-diagonals.
A prototype is
 $M= \Big| h_{i+j} \Big|_{i,j=1\ldots n}$,
with indeterminates $h_i$ in a commutative ring.

With one more set of indeterminates $a_i$, 
and an integer $k\in\Z$,
one defines the \emph{Hankel Pfaffian} $\Pf(\a,\h,n,k)$  to
be the Pfaffian of the antisymmetric matrix  
$\M(\a,\h,n,k)$ of order $2n$  
with entries $(a_i-a_j)h_{i+j-3+k}$ .  This is  the Pfaffian that
we shall study in this text. 
Such Pfaffians with $a_i=i$ or $a_i=q^i$ and special $h_i$
have been considered by Ishikawa, Tagawa, Zeng \cite{Ishi}.

Hankel matrices,
when the $h_i$ are identified with complete functions
of an alphabet of cardinality $n$, are related to
resultants, Bezoutians, orthogonal polynomials, 
continued fractions, etc \cite{CBMS}.
We show similarly in section 2 and section 5
that Hankel Pfaffians in complete functions
allow to express resultants, Bezoutians, $q$-discriminants,
and give several determinantal expressions of such Pfaffians.

The Hankel Pfaffian  $\Pf(\a,\h,n,k)$ can be studied by mere algebraic manipulations,
this is what we do in section 2. However, it is much more fruitful to use the action of
the symmetric group on the indeterminates $a_i$.  In \cite{LaPfaff}, we have shown 
how to diagonalize Pfaffians using Young's idempotents. 
In the present case,  it is more convenient to use 
the bases  of Kazhdan and Lusztig\cite{KL}. 
Theorem \ref{th:PfaffKL}
shows, indeed,  that $\Pf(\a,\h,n,k)$ is diagonal
in a pair of adjoint Kazhdan-Lusztig bases. 

Apart from the theory of symmetric functions,
we shall need  properties of  representations of 
the symmetric group, that we recall in section 3. Since the combinatorics of 
 Kazhdan and Lusztig bases are not well known, we give in this section more properties
 than  is needed  proper for the computation of Pfaffians.

\section{Symmetric functions}

We recall some properties of symmetric functions, 
following the conventions of \cite{CBMS}
rather than more classical ones as found in the book
of  I.G. Macdonald.

\subsection{Schur functions}

Given a sequence ${h_0=1,h_1,\ldots}$,
given an integer $n$ and $u,v \in\N^n$, one defines
the \emph{Schur function} $S_v$ to be the determinant
of $\Bigl[ h_{v_j+j-i} \Bigr]_{i,j=1\dots n}$,
and the \emph{skew Schur function} $S_{v/u}$ 
to be the determinant
of $\Bigl[ h_{v_j+j-i-u_i} \Bigr]_{i,j=1\dots n}$,
putting $h_i=0$ for $i<0$ (but this convention 
will be changed later).

Given two finite alphabets $\x=\{x_1,\dots,x_n\}$ and 
$\y=\{y_1,\dots,y_m\}$, the  
\emph{complete functions} $h_i(\x\moins \y)$
 of $\x\moins \y$ are defined by the generating series
$$ \prod_{i=1\dots m}(1-z y_i) 
   \prod_{i=1\dots n}(1-z x_i)^{-1} = \sum_0^\infty
   z^i h_i(\x\moins \y) \, . $$ 

Determinants of order $n$ in the complete functions of
$\x$ satisfy  \cite[Th.1.8.3]{CBMS}
\begin{equation}     \label{ProdSchur}
\det\bigl( h_{v_j+j-i+u_{n-i+}}(\x)   \bigr)
  = S_v(\x) S_u(\x) \ , \ u,v\in\N^n \, .
\end{equation}

Similarly, for  $r\ge 0$, $u,v\in\N^n$ such that 
$u\le r^n$, one has
\begin{equation}     \label{FactorSkewSchur}
S_{(v+r^n)/u}(\x) = S_v(\x)\,
  S_{r-u_n,\dots, r-u_1}(\x) \, .
\end{equation}

Schur functions of a difference of alphabets factorize,
when the components of $v$ are big enough \cite[Prop.1.4.3]{CBMS} 
\begin{equation}     \label{FactorSchur}
 S_{v+m^n}(\x-\y) = \prod_{i=1\dots n}\prod_{j=1..m}
    (x_i\moins y_j) \, S_v(\x)  \ ,\ v\in\N^n\, .
\end{equation}

\subsection{Invariance by translation of indices}

Given $\x$ of cardinality $n$,
the sequence $h_i(\x)$ is a recurrent sequence
\begin{equation}   \label{Recur_h}
   \sum_{i=0}^n  (\moins 1)^i e_i(\x)\, h_{k-i}(\x) = 0
 \qquad \text{ for   } k\geq n
\end{equation}
that one can extend, following Wronski,
 into a recurrent sequence  $\{ h_k(\x):\,
k\in\Z\}$ by requiring relation  \ref{Recur_h}  for all $k\in\Z$
and imposing the initial conditions
$h_{-1}(\x)=0=\cdots = h_{1-n}(\x)$ \cite{HouMu}.

From now on, the notation $h_k(\x)$,
as well as the different determinants in the $h_k(\x)$,
  will use this convention.
For example, a \emph{skew Schur
function} is defined for any pair $v,u\in\Z^n$~:
$S_{v/u}(\x) = \det( h_{v_j-u_i+j-i}(\x) )$.
In fact, one has
$$   \frac{1}{x_1\dots x_n} S_{v/u}(\x)  =
    S_{ (v-1^n)/u}(\x) = S_{v/(u+1^n)}(\x)  \, ,$$
so that, up to powers of $x_1\dots x_n$, one can recover 
indices in $\N^n$.   

The properties of determinants of order $n$ in the $h_i(\x)$
extend without further ado. For any $u,v\in\Z^n$, any $r\in\Z^n$, one has
\begin{equation}     \label{FactorSchur2}
S_{(v)/u}(\x) = S_{v+r^n}(\x)\,
  S_{r-u_n,\dots, r-u_1}(\x) \, .
\end{equation}
For example, for $n=3$, one has $h_{-1}(\x)=0$,
  $h_{-2}(\x)=0$, $h_{-3}(\x)= (x_1x_2x_3)^{-1}$ and 
$$ S_{023/001}(\x) =
\begin{vmatrix} h_0(\x) & h_3(\x) & h_5(\x) \\
               h_{-1}(\x) & h_2(\x) & h_4(\x) \\
            h_{-3}(\x) & h_0(\x) & h_2(\x)
\end{vmatrix} $$
factorizes into 
$ \left( x_1^{-1}\plus x_2^{-1}\plus x_3^{-1}\right) 
 S_{023}(\x)$,
but this is not the case of the determinant
$\begin{vmatrix} h_0(\x) & h_3(\x) & h_5(\x) \\
               0 & h_2(\x) & h_4(\x) \\
               0 & h_0(\x) & h_2(\x)
\end{vmatrix} $ corresponding to the conventions
$h_k(\x)=0$ for $k<0$.

More generally, given any alphabet $\y=\{y_1,\ldots, y_m\}$,
then $h_i(\x-\y)$ is a recursive sequence satisfying the same
recursion (\ref{Recur_h}), and therefore can be extended
to negative indices. The corresponding skew Schur functions
still satisfy, for $u,v\in\Z^n$,
\begin{equation}
  \frac{1}{x_1\dots x_n} S_{v/u}(\x)  =
    S_{ (v-1^n)/u}(\x) = S_{v/(u+1^n)}(\x)  \, .
\end{equation}

For example, for $n=2=m$, one has
$h_1(\x\moins\y)= x_1+x_2-y_1-y_2$, 
$h_0(\x\moins\y)= 1-y_1y_2(x_1x_2)^{-1}$,
$h_{-1}(\x\moins\y)= (y_1\plus y_2) (x_1x_2)^{-1}
  -y_1y_2(x_1\plus x_2) (x_1x_2)^{-2}$
and
$$ S_{02/00}(\x\moins \y) =
 \begin{vmatrix}  h_0(\x\moins\y) & h_3(\x\moins\y) \\
       h_{-1}(\x\moins\y) & h_2(\x\moins\y)
\end{vmatrix}
= \begin{vmatrix} 1 -e_2^y e_2^{-1} & h_3- e_1^y h_2+e_2^y h_1 \\
    e_1^y e_2^{-1} -e_2^y e_1e_2^{-2} & h_2- e_1^y h_1+e_2^y 
\end{vmatrix}  $$ 
$$ = (x_1x_2)^{-2} S_{24/00}(\x\moins \y) = (x_1x_2)^{-2}
  \, R(\x,\y) (x_1^2\plus x_1x_2\plus x_2^2) \, ,$$
writing $e_i,h_i$ for the functions of $\x$, and 
$e_1^y,e_2^y$ for those of $\y$.

\subsection{Bezoutians and resultants}
Given $\x=\{x_1,\dots,x_n\}$, the remainder of a polynomial
$f(y)$ modulo $S_n(y\moins \x)$ is the only polynomial
$\cR_\x f$ of degree $\le n\moins 1$ such that
\begin{equation}     \label{Reste}
 \cR_\x f( x_i) = f(x_i)   \ , \quad i=1,\ldots,n \, .
\end{equation}
Thus, the definition of the remainder can be extended to
any function $f(y)$, in particular, can be extended 
\cite{Boliya} to polynomials in $y,y^{-1}$, 
by requiring relations \ref{Reste}.

Similarly, the \emph{Bezoutian}  $\Bez_\x(f)$ of
a function $f(y,z)$ is
the matrix  
$\Bez_\x(f) =\Bigl[  b_{ij} \Bigr]_{i,j=0\ldots n-1}$,
where $\sum_{i,j=0}^n b_{ij} z^{n-1-i} y^{n-j-1}$
is the remainder of
$f(y,z)$ modulo $S_n(z-\x)$ and modulo $S_n(y-\x)$
\cite[Th.3.4.1]{CBMS,LaPragacz}.

Given another alphabet $\mathbf c=\{c_1,\ldots,c_m\}$,
the \emph{resultant} $R(\x,\mathbf c)$ is defined to be
$$ R(\x,\mathbf c) =\prod_{i=1}^n \prod_{j=1}^m (x_i-c_j) \, .$$
It is also \cite[Th.3.2.1]{CBMS} equal  to the Schur function
$S_{m^n}(\x-\mathbf c)$.

\begin{lemma}
The determinant of
$\Bez_\x\Bigl(  S_{n-1}(y+z) S_m(z\moins \mathbf c) \Bigr) $
is equal to the resultant $R(\x,\mathbf c)$.
\end{lemma}
\proof  Instead of expanding the double remainder in the basis
$y^jz^i$, let us choose the basis $y^j S_i(z\moins \mathbf c)$,
$i,j=0\ldots n\moins 1$.
One has
\begin{eqnarray*}
S_{n-1}(y+z) S_m(z\moins \mathbf c)
&=&
\sum_{j=0}^{n-1} y^{n-1-j} z^j S_m(z\moins \mathbf c) =
\sum_{j=0}^{n-1} y^{n-1-j} S_{m+j}(z\moins \mathbf c) \\
&=&
\sum_{j=0}^{n-1} y^{n-1-j} S_{m+j}\Bigl( ( z\moins\x)
             + (\x\moins \mathbf c) \Bigr)      \\
&\equiv& \sum_{j=0}^{n-1} \sum_{i=0}^{n-1}
  y^{n-1-j} S_i(z\moins\x) S_{m+j-i}(\x\moins \mathbf c) \, .
\end{eqnarray*}
Hence, the Bezoutian, expressed in these bases, is the matrix
$\Bigl[ S_{m+j-i}(\x\moins \mathbf c)  \Bigr]$,
the determinant of which is equal to $S_{m^n}(\x-\mathbf c)$.  \QED

Using $\{c_1,\ldots,c_m,0\}$ instead of $\mathbf c$,
one obtains that the determinant of
$\Bez_\x\Bigl(  S_{n-1}(y+z)\, z S_m(z\moins \mathbf c) \Bigr) $
is equal to $x_1\dots x_n R(\x,\mathbf c) =
S_{(m+1)^n}(\x-\mathbf c)$, and more generally, that 
\begin{equation}   \label{Bez1}
\det\left( \Bez_\x\Bigl(  S_{n-1}(y+z)\, z^k 
     S_m(z\moins \mathbf c) \Bigr)  \right)
 =S_{(m+k)^n}(\x-\mathbf c) \, ,
\end{equation}
the equality being valid for $k\in\Z$, 
once its is checked for a single value of $k$.

For example, for $n=2=m$, the matrices of the remainders
in the basis $\{1,y\} \otimes \{1,z\moins \x\}$ are,
for $k=0,\, \moins1,\, \moins 2$ respectively,
$$ \begin{bmatrix}
  S_1(\x \moins \mathbf c) & S_2(\x \moins \mathbf c)\\
  S_2(\x \moins \mathbf c) & S_3(\x \moins \mathbf c)
\end{bmatrix} \, ,\
\begin{bmatrix}  1-\frac{c_1c_2}{x_1x_2} & S_1(\x \moins \mathbf c)
 \\ S_1(\x \moins \mathbf c) & S_2(\x \moins \mathbf c)
\end{bmatrix} \, ,\
\begin{bmatrix}
 \frac{y_1+y_2}{x_1x_2} -\frac{y_1y_2(x_1\plus x_2)}{(x_1x_2)^2}
          &1-\frac{c_1c_2}{x_1x_2} \\
1-\frac{c_1c_2}{x_1x_2} & S_1(\x \moins \mathbf c)
\end{bmatrix} \, .
$$

\subsection{Discriminants}

Let $f(y)=R(y,\x)$. Write the derivative $f'(y)$ of $f(y)$
in the factorized form $f'(y)= n R(y,\x^{der})$,
using the alphabet $\x^{der}$ of roots of 
$f'(y)$. The logarithmic derivative of $f'(y)$
shows that $p_k(\x)= n h_k(\x\moins \x^{der})$ 
for any $k\ge 0$. 
Since for any $i$, $\{ x_i^r\}$ is 
a recursive sequence satisfying
(\ref{Recur_h}), $\{ p_r(\x)= x_1^r\plus \dots\plus
x_n^r ,\
  r\in\Z\}$ is a recursive sequence satisfying (\ref{Recur_h}).
Thus the equality $p_k(\x)= h_k(\x\moins \x^{der})$
can be extended to any $k\in\Z$ by taking the preceding
conventions for $h_k(\x)$ and $p_k(\x)$, $k\in\Z$.

The resultant of $R(\x,\x^der)$ is equal to 
$S_{(n-1)^n}(\x\moins \x^{der})$, that is, 
to the determinant
$ \Bigl| \frac{\s 1}{\s n}p_{n-1+j-i}(\x) \Bigr|_{i,j=1\dots n}$, which is equal to 
$(\moins 1)^{\binom{n}{2}} \cD(\x,1)$,
where $\cD(\x,1)$ is the square of the Vandermonde in $\x$,
called the \emph{discriminant}. 

More generally, the resultant $R(\x, q\x)=
  \prod_{i,j=1\dots n} (x_i -qx_j) =S_{n^n}(\x-q\x)$ 
is equal to 
$(\moins 1)^{\binom{n}{2}} (1\moins q)^n x_1\dots x_n$
times the \emph{$q$-discriminant} 
$$ \cD(x,q)= \prod_{1\le i<j\le n} (x_i-qx_j) (x_j-qx_i)\, .$$

\section{Determinantal expressions}

It is clear that $\Pf(a,h,n,0)$ is of degree $n$ in the 
variables $h_0,h_1,\ldots$, and thus expands as a sum
of Schur functions of index in $\N^n$. One can therefore
introduce $\x=\{x_1,\dots,x_n\}$ and specialize
each $h_i$ to the complete function $h_i(\x)$ without loss
of information (thanks to homogeneity, $h_0(\x)=1$ creates
no problem). Moreover, we have seen in the preceding section
that shifting the indices $h_i(\x)\to h_{i+k}(\x)$
with a fixed $k\in\Z$ multiplies the skew Schur functions
of $\x$ by a factor $(x_1\dots x_n$. Hence, it is
easy to pass from the Pfaffian in 
$(a_i\moins a_j) h_{i+j-3}(\x)$ to the Pfaffian in
$(a_i\moins a_j) h_{i+j-3+k}(\x)$.

Let $E(\x)$ be the matrix of order $2n$ with entries
the signed \emph{elementary symmetric functions} 
$  (\moins 1)^{j-i}e_{j-i}(\x)$ (defined to be $0$ 
for negative indices). 
It is straightforward that sums of the type 
$\sum_{i=0}^j (\moins 1)^i e_{j-i}(\x) h_{k+i}(\x)$,
for $k\ge n\moins 1$, and $j\ge 0$,  
are equal to \emph{hook Schur functions} $S_{1^j,k}(\x)$.

Using this property, one checks the following proposition
by decomposing linearly the Pfaffian matrix
according to the $a_i$'s.

\begin{proposition}
  The matrix $ E^{tr}(\x)\, \M(\a,h(\x),n)\, E(\x)$
  is such that its submatrix on rows and columns 
  $n\plus 1,\ldots,2n$ is null.
 The submatrix on rows $1,\ldots,n$ and columns 
  $n\plus 1,\ldots,2n$, denoted $\Coin(\a,\x,n)$,
 is equal to
 \begin{equation}   \label{Filtre Coin}
  \Coin(\a,\x,n)
   = \sum_{r=1}^{2n} a_r \, \Big[ (\moins 1)^{n+r+i+j}
     e_{n-r+j}(\x)\, S_{1^{i-1},\, r-2}(\x)  \Big] \, .
\end{equation}
\end{proposition}

For example, for $n=2$, one has
\begin{multline*}
\begin{bmatrix} 1 & 0&0&0\\ -e_1 &1& 0&0\\ e_2 &-e_1 &1& 0\\
  0 &e_2 &-e_1 &1 \\  \end{bmatrix}
\begin{bmatrix}
0& (a_1-a_2)S_{0}& (a_1-a_3)S_{1}& (a_1-a_4)S_{
2}\\ (a_2-a_1)S_{0}& 0& (a_2-a_3)S_{2}& (a_2-a_4)S_{
3}\\ (a_3-a_1)S_{1}& (a_3-a_2)S_{2}& 0& (a_3-a_4)S_{
4}\\ (a_4-a_1)S_{2}& (a_4-a_2)S_{3}& (a_4-a_3)S_{4
}& 0
\end{bmatrix} \\
\times \begin{bmatrix} 1 & -e_1 & e_2 & 0\\  0 & 1 & -e_1 & e_2\\
  0 & 0 &1 & -e_1 \\ 0&0&0&1 \end{bmatrix}  =\\ 
\begin{bmatrix}
0& (a_1\moins a_2) S_0& (a_2\moins a_3) S_1
& (a_3\moins a_4) S_2+( a_3\moins a_2) S_{11} \\ 
(a_2 \moins a_1)   S_0& 0& (a_3\moins a_1)
S_{11} & (a_4\moins  a_3) S_{12} \\
 (a_3\moins a_2) S_1
& (a_1 \moins a_3)  S_{11} & 0& 0\\
( a_4 \moins a_3) S_2 +(a_2 \moins a_3)
S_{11} & (a_3 \moins a_4)  S_{12} & 0& 0
\end{bmatrix} \, .    
\end{multline*}
The $2\times 2$ North-East corner $\Coin(\a,\x,2)$ expands as
\begin{multline*}
 a_1 \begin{bmatrix} \moins e_2 S_{0,-1} & 0\\ 
                  \moins  e_2 S_{1,-1} & 0\end{bmatrix} 
 + a_2 \begin{bmatrix} e_1 S_{0,0} & \moins e_2S_{0,0}\\
               e_1 S_{1,0} & \moins e_2S_{1,0} \end{bmatrix}
+ a_3 \begin{bmatrix} \moins e_0 S_{0,1} & e_1 S_{0,1}\\
     e_0 S_{1,1} & \moins e_1 S_{1,1}  \end{bmatrix}
+ a_4 \begin{bmatrix} 0 & \moins e_0 S_{0,2} \\
       0 & e_0 S_{1,2} \end{bmatrix} \\
=a_1 \begin{bmatrix} 0 &0\\ -e_2 & 0 \end{bmatrix}
 + a_2 \begin{bmatrix}  e_1 &\moins e_2\\ 0 & 0 \end{bmatrix}
 + a_3 \begin{bmatrix}  \moins e_1 & e_1 S_1
        \\ e_2 & \moins e_2 S_1 \end{bmatrix}
+ a_4 \begin{bmatrix}  0 &  \moins S_2 \\ 0 & 
                 e_1 S_2\end{bmatrix} \, .
\end{multline*}

One can also write the Pfaffian 
as the determinant
of a $2n\times 2n$ matrix in different manners, as shows the
next result.

\begin{proposition}
 For a given $n$, let
 $\M_h$ be the matrix of order $2n$ with $i$-th row
 $$ [ h_{n+1-i}(\x),\dots, h_{2n-i}(\x),
    a_{n-i+1} h_{n+1-i}(\x),\dots, a_{n-i+1} h_{2n-i}(\x)]\, ,$$
 let $\M_e$ be the matrix obtained by changing each $h_i(\x)$
 to $e_i(\x)$ in $\M_h$, 
 and finally, let $\M_x$ be the matrix with
 $i$-th row 
 $$ [ x_1^{2n-i},\dots, x_n^{2n-i}, a_{2n-i+1}x_1^{2n-i},
       \dots, a_{2n-i+1} x_n^{2n-i}]\, .$$
Then $\Pf(\a, h(\x), n, 2\moins n)$ is equal to
the determinant of $\M_h$, of $\M_e$,
and equal to the quotient of the determinant of $\M_x$
by $\Delta(\x)^2$.
\end{proposition}
\proof It results from the analysis in \cite{LaPfaff}
that  the Pfaffian is determined by  specializing
$\a$ to a permutation of $1^n 0^n$ (in fact, Yamanouchi
words suffice). In that case, the Pfaffian becomes 
equal, up to a sign, to some determinant of order of $n$ 
(denoted $g(i\dots j|k\dots \ell)$ in \cite{LaPfaff}).
In the present case, this determinant is 
of the type  $\Big| h_{u_i+v_j+i+j-n-1}(\x) \Bigr|$,
with $u,v$ two increasing partitions in  $\N^n$.
But such a determinant is equal to
$(\moins 1)^{\binom{n}{2}} S_u(\x) S_v(\x))$ thanks to
(\ref{ProdSchur}). 
Similarly, the specializations of $\M_h, \M_e$ or $\M_x$
give rise to products of two Schur functions.  \QED

For example, for $n=2$, 
the three matrices appearing in the proposition are
$$ 
\begin{bmatrix}   h_2  &  h_3 &   a_4 h_2 &   a_4 h_3\\
         h_1 &   h_2  &  a_3 h_1  &  a_3 h_2\\
         1    & h_1  &   a_2   &   a_2 h_1\\
        0    & 1    &   0     &  a_1
\end{bmatrix}  \ , \
\begin{bmatrix}  e_2 &   0 &    a_4 e_2  &    0\\
                e_1 &   e_2  &  a_3 e_1  &  a_3 e_2\\
                1   &  e_1  &   a_2  &    a_2 e_1\\
           0    & 1    &   0     &  a_1
\end{bmatrix} \ , \
 \begin{bmatrix} x_1^3&   x_2^3 &  a_4 x_1^3& a_4 x_2^3 \\
               x_1^2&   x_2^2 &  a_3 x_1^2& a_3 x_2^2 \\ 
              x_1&   x_2 &  a_2 x_1 & a_4 x_2 \\
              1 & 1  & a_1 & a_1 
\end{bmatrix}   \, ,
$$ 
and the determinant of the first two matrices is equal to
$$ \Pf(\a,h(\x) ,2,0)=
  (a_3 - a_2) (a_4 - a_1) S_{22}(\x)
  - (a_2 - a_1) (a_4 - a_3) S_{13}(\x)  \, .
$$

One can extend the preceding property by shifting
indices: $h_i\to h_{i+r}$, and considering functions
of $\x-\y$ instead of $\x$. 
Thus given $\y=\{y_1,\dots,y_m\}$ of cardinality $m$, 
and any $r\in\Z$
let $\M_h(\a,\x-\y,r)$ be the matrix of order $2n$
with $i$-th row
\begin{multline*}
 [ h_{n+1-i+r}(\x-\y),\dots, h_{2n-i+r}(\x-\y),  \\
    a_{n-i+1} h_{n+1-i+r}(\x-\y),\dots, 
          a_{n-i+1} h_{2n-i+r}(\x-\y)]\, . 
\end{multline*}
Thanks to (\ref{FactorSchur2}), each minor on the first
$n$ columns, or last $n$ columns of $\M_h(\a,\x-\y,r)$
is equal to the product of the same minor 
of $\M_h(\a,\x-\y,0)$  by $R(\x,\y)\, (x_1\dots x_n)^{r-m}$.
Similarly, 
$$ \Pf(\a,h(\x-\y),k)=  \Pf(\a,h(\x),0) 
         R(\x,\y)\, (x_1\dots x_n)^{k-m}  \, ,$$
since the Pfaffian is a linear combination of 
Schur functions.
Hence the preceding proposition entails~:
\begin{theorem}
Given two finite alphabets $\x$ of cardinality $n$,
$\y$ of cardinality $m$, and two integers $k,r\in\Z$
such that $k+m+n-2r-2=0$,
then one has
\begin{eqnarray}
\Pf(\a, h(\x-\y),k)  &=&  (x_1\dots x_n)^{k-m} \, R(\x,\y)\,
                          \Pf(\a, h(\x), 0)  \\
 \label{PfaffByMatAleph}
  &=& \frac{1}{R(\x,\y)}\, \det\bigl( \M_h(\a,\x-\y,r)\bigr) \, .
\end{eqnarray}
\end{theorem}

For example, for $n=2$, $m=1$, the matrix
$$   \M_h(\a,\x\moins \y,0)= 
\begin{bmatrix}   h_2(\x\moins\y)  &  h_3(\x\moins\y)
            &   a_4 h_2(\x\moins\y) &   a_4 h_3(\x\moins\y)\\
         h_1(\x\moins\y) &   h_2(\x\moins\y)  &  
   a_3 h_1(\x\moins\y)  &  a_3 h_2(\x\moins\y)\\
     1    & h_1(\x\moins\y)  &   a_2   &   a_2 h_1(\x\moins\y)\\
    y_1(x_1x_2)^{-1}    & 1    &  a_1y_1(x_1x_2)^{-1}     &  a_1
\end{bmatrix}  $$
has determinant equal to
$$  (x_1-y_1)^2 (x_2-y_1)^2 
   \left((a_1-a_3) (a_2-a_4)  -
    (a_1-a_2) (a_3-a_4) \frac{(x_1+x_2)^2}{x_1x_2}  \right)\, .$$

There is another way to evaluate a Hankel Pfaffian,
when $\a$ is specialized to $\q=\{1,q,q^2,\dots, q^{2n-1}\}$.
The next theorem shows that in that case the matrix
$\Coin\Bigl(\a,\x,n)$ coincides with a Bezoutian.

\begin{theorem}
 Given $\x$ of cardinality $n$,  one has for $k\in \Z$
\begin{equation}   \label{CompareCoinBezout}
  \Coin\Bigl( (q^{i-1}-q^{j-1}) h_{i+j+k-n-2}(\x)  \Bigr)
 = \Bez_\x\biggl(- z^k S_{n-1}(qz+y) S_n(qz-\x)  \biggr) \, .
\end{equation}
\end{theorem}
\proof
 Suppose $k\ge 0$ and expand
\begin{eqnarray*}
 \moins z^k S_{n-1}(qz+y) S_n(qz-\x) &=&
  -\sum_{j=0}^{n-1} y^j z^{k+n-1-j} q^{n-1-j} S_n(qz-\x) \\
 &=&
 -\sum_{j=0}^{n-1} q^{2n-1-j} y^j S_{2n+k-j-1}(z-q^{-1}\x)\, .
\end{eqnarray*}
Using the expression of the remainder as a hook Schur function
given in \cite[Th. 3.2.1]{CBMS}, one has
\begin{equation*}
  S_m(z-B) \equiv 
   \sum_{j=0}^{m} \sum_{i=0}^{n-1}
   (\moins z)^i S_{1^{n-1-i},m-n+1-j}(\x) 
       S_j(\moins q^{-1}\x)   \ \mod R(z,\x)
\end{equation*}
and one obtains that the Bezoutian is a matrix
with entries equal to hook Schur functions of $\x$
times functions
$S_j(\moins q^{-1}\x) = (\moins 1)^j q^{-j} e_j(\x)$.
More precisely, filtering the Bezoutian according
to powers of $q$, 
one recognizes in this filtration exactly the
filtration of $\Coin(\a,n,k)$ according to $a_1,\dots,a_{2n}$
obtained from (\ref{Filtre Coin}).
The expression remains valid for $k<0$ because the
entries of the Bezoutian for variable $k$ form a recursive
sequence with the same characteristic polynomial as
the sequence $h_k(\x)$, $k\in\Z$. \QED

For example, for $n=2$ and $k=0$, one has
\begin{multline*}
 \Bez_\x\bigl(-(qz+y) S_2(qz-\x) \bigr)=
 \Bez_\x\bigl( -S_3(qz-\x)-yS_2(qz-\x) \bigr)  \\
 = \begin{bmatrix} 0 &0\\ -e_2 & 0 \end{bmatrix}
 + q \begin{bmatrix}  e_1 &\moins e_2\\ 0 & 0 \end{bmatrix}
 + q^2 \begin{bmatrix}  \moins e_1 & e_1 S_1
        \\ e_2 & \moins e_2 S_1 \end{bmatrix}
+ q^3 \begin{bmatrix}  0 &  \moins S_2 \\ 0 &
                 e_1 S_2\end{bmatrix} \, .
\end{multline*}

Since $-q S_{n-1}(qz+y) S_n(qz-\x)
   = -S_{n-1}(z+q^{-1}y)S_n(z-q^{-1}\x)$, 
the determinant of
$\Bez_\x\biggl( S_{n-1}(qz+y) S_n(qz-\x)  \biggr)$
is equal to the resultant $R(\x,q^{-1}\x)$ up to a power of $q$
and a sign. 

Controlling the power of $q$, and using the invariance of
$$ \det\biggl(\Bez_\x\bigl( -z^k S_{n-1}(qz+y) S_n(qz-\x)
          \bigr) \biggr) (x_1\dots x_n)^{-k} $$ 
with respect to $k$,
one obtains the following property.

\begin{theorem}     \label{qDiscri}
 Given $\x$ of cardinality $n$,  one has for $k\in \Z$
\begin{eqnarray}
  \Pf\Bigl( (q^{i-1}\moins q^{j-1}) h_{i+j+k-n-1}(\x)  \Bigr)
 &=& (\moins 1)^{\binom{n}{2}} S_{n^n}((1\moins q)\x) (x_1\dots x_n)^k \\
 &=& (\moins 1)^{\binom{n}{2}} (1\moins q)^n
    (x_1\dots x_n)^{k+1} \cD_\x(q)   .
\end{eqnarray}
\end{theorem}

For example, for  $n=3$, and $k=-1$, one has
$h_{-1}(\x)=0= h_{-2}(\x)$ and the determinant of
the skew-symmetric matrix 
$$  \begin{bmatrix}
 0& 0& 0& 1-q^{3}& (1-q^{4})S_{1}& (1-q^{5})S_{2}\\
0& 0& q-q^{2}& (q-q^{3})S_{1}& (q-q^{4})S_{2}& (q-q^{5})S_{3
}\\ 0& q^{2}-q& 0& (q^{2}-q^{3})S_{2}& (q^{2}-q^{4})S_{3
}& (q^{2}-q^{5})S_{4}\\
 q^{3}-1& (q^{3}-q)S_{1}& (q^{3}-
q^{2})S_{2}& 0& (q^{3}-q^{4})S_{4}& (q^{3}-q^{5})S_{5}\\
(q^{4}-1)S_{1}& (q^{4}-q)S_{2}& (q^{4}-q^{2})S_{3}& (q^{4}-q
^{3})S_{4}& 0& (q^{4}-q^{5})S_{6}\\
 (q^{5}-1)S_{2}& (q^{
5}-q)S_{3}& (q^{5}-q^{2})S_{4}& (q^{5}-q^{3})S_{5}& (q^{5}-q
^{4})S_{6}& 0
\end{bmatrix}
$$
is equal to the square of $ (1\moins q)^3 q^3\, \cD_\x(q)$.

\section{Representations of the symmetric group}

To understand the dependency in $\a$ of the Pfaffian 
$\Pf(\a,\h,n)$, we need to use the theory of representations.
Irreducible representations of the symmetric group 
$\mfS_n$ over $\C$ are in bijection with partitions of $n$.
One usually indexes bases by \emph{standard Young tableaux}
of a given shape $\l$. 
The tableaux of shape $\l$ can be considered as the vertices of
a graph, two tableaux being connected 
by an edge of label $s_i$ if the two tableaux 
differ by the transposition of $i,i\plus 1$. 

Interpreting tableaux of shape $\l$ 
as products of Vandermonde  determinants,
each column $u=[u_1,\dots,r_r]$ giving rise to the 
Vandermonde 
$\Delta^x(u)= \prod_{1\leq i<j\leq r} (x_{u_i} -x_{u_j})$,
one obtains the \emph{Specht basis} of the irreducible 
representation of index $\l$ of the symmetric group.
More generally, we shall call Specht basis
any image of this basis in another copy of the same
representation. 

In \cite{LaPfaff}, we have used a \emph{Young basis} rather
than a Specht basis to expand a Pfaffian of the type
$\Pf( (a_i-a_j) g_{i,j})$, $1\leq i<j\leq 2n$, $g_{i,j}=g_{j,i}$,
observing that three symmetric groups are involved:
the symmetric group permuting the $a_i$, the symmetric group
acting on $g_{i,j}$, and the \emph{diagonal} group
acting simultaneously on the indices of $a_i$ and $g_{i,j}$.

In the case of a Hankel Pfaffian, it will be more illuminating
to use several \emph{Kazhdan-Lusztig bases}, corresponding to 
different spaces of polynomials.  The original
constructions of Kazhdan and Lusztig stand at the level of
the Hecke algebra. Unfortunately,  general irreducible representations 
are still not fully explicit. However the case of interest for
Pfaffians is the case corresponding to Gra{\ss}mannians 
\cite{Torun}, that is, the case of rectangular 
partitions of the type $[n,n]$, or $[2^n]$ 
that one can find in the literature under many disguises.

We shall need only a pair of bases,
but prefer to be more complete and describe 
the Kazhdan-Lusztig bases of some other realizations of the same representation.

\subsection{Combinatorial objects}

Bases of irreducible representations of the symmetric group are
usually encoded by \emph{standard Young tableaux} of a given shape.
In our case, the shape will be $[n,n]$ or its transpose $[2^n]$.

From a $2\times n$ Young tableau, one reads two partitions,
by subtracting to  the bottom row, as a vector, the vector
  $[1,2,\dots, n\moins 1]$, and by subtracting
to $[n\plus 1,\dots, 2n]$ the top row.
These two partitions are contained in the staircase partition
$[n\moins 1,\dots,1,0]$. We shall label bases by the partition $\l$
(written decreasingly) corresponding to the bottom row of the tableau.

Thus the Young tableau
$$\young{ 3& 6& 7& 9 & 11 &12\cr
  1 &2 &4 &5 &8 & 10 \cr}  \ 
 \Rightarrow 
\begin{array}{rl} {} [7,\dots,12]-[3,6,7,9,11,12] &= [4,2,2,1,0,0] \\
      {} [1,2,4,5,8,10]-[1,\dots,6] &= [0,0,1,1,3,4] 
\end{array}
$$
will be replaced by $\l= [4,3,1,1]$.

To $\l$ one also associate a \emph{skew partition} 
 $\boxed{\s\l} = \bigl( [(n\moins 1)^n] + 
    \l\omega \bigr)\, /\l^\sim$, where $\l\omega$ means
the increasing reordering of $\l$.
For the running example, it is 
$$ \boxed{\l}=  \bigl([5^6]+ [0,0,1,1,3,4]\bigr)/[0,0,1,2,2,4]
  = [5,5,6,6,8,9]/[0,0,1,2,2,4]\, .  $$

Reading the border of the diagram of $\l$, one obtains
a \emph{Yamanouchi word} that one can represent planarly 
as a \emph{Dyck path}, $1$ standing for a North-East step,
$0$ a South-East step.
For $n=4$, $\l=[3,1]$ (figured in red) one has

$$  
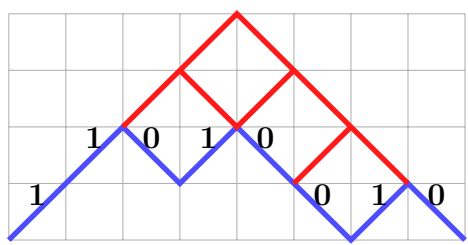
\begin{tikzpicture}[scale=0.75]
\draw[color=black!30] (0,0) grid (8,4);

\draw[Pas] (0,0)--(1,1) --(2,2)--(3,1)--(4,2)--(5,1)--(6,0)--(7,1)-- (8,0);
\draw[Pas2] (2,2)--(3,3)--(4,4)--(5,3)--(6,2)--(7,1);
\draw[Pas2] (3,3)--(4,2)--(5,3);
\draw[Pas2] (5,1)--(6,2);

    \node at (0.5,0.8) {$\mathbf 1$};
    \node at (1.5,1.8) {$\mathbf 1$};
    \node at (2.5,1.8) {$\mathbf 0$};
    \node at (3.5,1.8) {$\mathbf 1$};
    \node at (4.5,1.8) {$\mathbf 0$};
    \node at (5.5,0.8) {$\mathbf 0$};
    \node at (6.5,0.8) {$\mathbf 1$};
    \node at (7.5,0.8) {$\mathbf 0$};
\end{tikzpicture}
  \quad \raisebox{40pt}{$ \begin{array}{c}
    \lambda =[3,1]  \\[4pt]
  {Yamanouchi}\  \mathbf{[1,1,0,1,0,0,1,0]} \end{array} $}
 $$

Pairing successively in the Yamanouchi word 
$1\dots 0$ treated as opening and
closing parentheses, one obtains a \emph{link pattern}.
To a link between positions $i$ and $j$ one associates
a factor $(a_i \moins a_j)$. Let 
$\varphi^a(\l)$ be the product of all such factors
for the link pattern associated to $\l$.
Equivalently, one labels the steps of the path by $1,2,\dots, 2n$,
each factor $(a_i \moins a_j)$ corresponding to  paired steps. 
$$ 
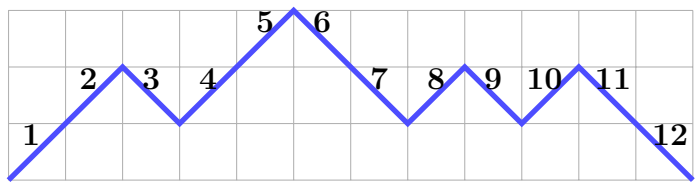
\begin{tikzpicture}[scale=0.75]
    \draw[color=black!30] (0,0) grid (12,3);
  \draw[Pas] (0,0)--(1,1) --(2,2)--(3,1)--(4,2)--(5,3)--(6,2)--(7,1)
--(8,2)--(9,1)--(10,2)--(11,1)--(12,0);

    \node at (0.4,0.8) {$\mathbf 1$};
    \node at (1.4,1.8) {$\mathbf 2$};
    \node at (2.5,1.8) {$\mathbf 3$};
    \node at (3.5,1.8) {$\mathbf 4$};
     \node at (4.5,2.8) {$\mathbf 5$};
     \node at (5.5,2.8) {$\mathbf 6$};
     \node at (6.5,1.8) {$\mathbf 7$};
     \node at (7.5,1.8) {$\mathbf 8$};
     \node at (8.5,1.8) {$\mathbf 9$};
     \node at (9.4,1.8) {$\mathbf{10}$};
     \node at (10.6,1.8) {$\mathbf{11}$};
     \node at (11.6,0.8) {$\mathbf{12}$};
\end{tikzpicture}$$
$$ \varphi^a([4,3,1,1]) =
(a_1\moins a_{12})(a_2\moins a_3)(a_4\moins a_7)(a_5\moins a_6)(a_8\moins a_9)
(a_{10}\moins a_{11})    \, $$

Let $\psi(\l)$ be the vector obtained by labeling 
$0,2,\ldots,2n\moins 2$ the successive increasing steps of
the Dyck path, and 
labeling each descending step by the label
of the step to which it is paired. 
$$ 
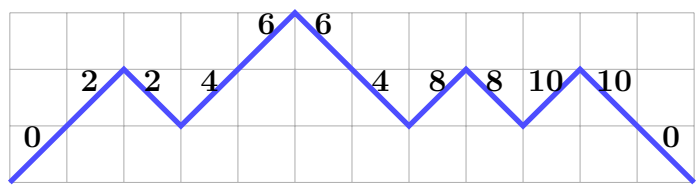
\begin{tikzpicture}[scale=0.75]
    \draw[color=black!30] (0,0) grid (12,3);
  \draw[Pas] (0,0)--(1,1) --(2,2)--(3,1)--(4,2)--(5,3)--(6,2)--(7,1)--(8,2)--(9,1)--(10,2)--(11,1)--(12,0);

    \node at (0.4,0.8) {$\mathbf 0$};
    \node at (1.4,1.8) {$\mathbf 2$};
    \node at (2.5,1.8) {$\mathbf 2$};
    \node at (3.5,1.8) {$\mathbf 4$};
     \node at (4.5,2.8) {$\mathbf 6$};
     \node at (5.5,2.8) {$\mathbf 6$};
     \node at (6.5,1.8) {$\mathbf 4$};
     \node at (7.5,1.8) {$\mathbf 8$};
     \node at (8.5,1.8) {$\mathbf 8$};
     \node at (9.4,1.8) {$\mathbf{10} $};
     \node at (10.6,1.8) {$\mathbf{10} $};
     \node at (11.6,0.8) {$\mathbf 0$};
\end{tikzpicture}$$
$$ \psi([4,3,1,1])= [0, 2, 2, 4, 6, 6, 4, 8, 8, 10, 10, 0] \, . $$

Given $n$ and a partition $\l$ one labels the boxes
of the diagram of $\l$ by a pair of numbers.
The first one increases by $1$ when moving horizontally
rightwards, and decreases by $1$ when moving 
vertically downwards, starting from $n$ in the first box.
The second number is $0$ for the boxes in the corners,
$1$ for the new corners obtained by erasing the preceding
corners, and so on. Let us denote this bi-labelled diagram
 $\cD_\l$.
A similar construction is given in type $B$ by 
\cite{deGierPyatov}.
$$  \cD_{4311}=   \bigyoung{6,3 & 7,2 & 8,1 & 9,0 \cr
               5,2 & 6,1 & 7,0 \cr
               4,1 \cr
               3,0 \cr }  $$

The weights $\varphi^a(\l)$ and $\psi(\l)$ 
will be interpreted as dual Kazhdan-Lusztig bases,
while the diagrams $\cD_\l$ will be used to generate
several Kazhdan-Lusztig bases \cite{KiriLa}.

\subsection{Basis $\KL_\l^\Delta$}
One generates it from
$ \KL_0^\Delta:= \Delta^x(1\ldots n\, |\, n+1\ldots 2n)$.
The polynomial $ \KL_\l^\Delta$ is defined to be the image of
$\KL_0^\Delta$ under $\cD_\l$, the diagram being read by 
successive rows,
each entry $[i,k]$ being interpreted as  
$s_i - (1\plus k)^{-1}$.

For example, for $n=3$, one has  
\begin{eqnarray*}
\KL_{11}^\Delta&=&\Delta^x(123|456)(s_3\moins \demi)(s_2\moins 1)\\
&=& \Delta^x(123|456)\Bigl( s_3s_2 -s_3 -\demi(s_2\moins 1)\Bigr)\\
&=& \Delta^x(134|256) -\Delta^x(124|356) +\Delta^x(123|456) \\
&=& \Delta^x(234|156)   \, ,
\end{eqnarray*}
the last expression being due to the Pl\"ucker relations.

The full basis for $n=3$ is 

$$ \xymatrix@R=0.5cm@C=-1.8cm{
 &  &  & *{\GrTeXBox{
              \KL_0^\Delta=\Delta^x(123|456)
}}\arx2[d]& \\
 &  &  & *{\GrTeXBox{
          \KL_1^\Delta=\Delta^x(124|356) \moins\Delta^x(123|456)
}}\arx1[ld]\arx3[rd]& \\
 &  & *{\GrTeXBox{
     \KL_{11}^\Delta=\Delta^x(234|561)
}}\arx3[rd]&  & *{\GrTeXBox{
           \KL_2^\Delta=\Delta^x(345|612)
}}\arx1[ld]& \\
 &  &  & *{\GrTeXBox{
  \KL_{21}^\Delta=\Delta^x(235|461)\moins\Delta^x(234|561)
}}& \\
}$$

\subsection{Basis $\KL_\l^S$}      

The family of skew Schur functions 
$S_{\boxed{\s\l}}$, $\l\leq \rho$,
is the Specht basis of an irreducible representation 
of index $[2^n]$. One interprets now standard tableaux 
as skew Schur functions, instead of products of Vandermonde determinants,
keeping the same action of the symmetric group.

The Kazhdan-Lusztig basis $\KL_\l^S$  is obtained from 
$ \KL_0^S= S_{(n-1)^n}$  using the diagrams 
$\cD_\l$, $\l\leq \rho$, interpreting an entry $[i,k]$
as $s_i - (1\plus k)^{-1}$. 

For example, for $n=3$, the Specht basis is
$S_{\,\boxed{\s 0}}= S_{222}$, 
$S_{\,\boxed{\s 1}}= S_{223/001}$,
$S_{\,\boxed{\s 11}}= S_{233/002}$,
$S_{\,\boxed{\s 2}}= S_{224/011}$,
$S_{\,\boxed{\s 21}}= S_{234/012}$.

Consequently, 
\begin{eqnarray*}
 \KL_1^S &=& S_{222} (s_3-1)= S_{223/001}-S_{222}= S_{123} \\
 \KL_2^S &=& S_{222} (s_3-\demi)(s_4-1)= 
S_{224/011} - S_{223/001}+S_{222}= S_{114} \\
 \KL_{11}^S &=& S_{222} (s_3-\demi)(s_2-1)=
S_{233/002}- S_{223/001}+S_{222}= S_{033} \\
 \KL_{21}^S &=& S_{222} (s_3-\demi)(s_2-1)(s_4-1)=
S_{234/012} -S_{224/011} \\
  & & \hspace*{60pt}  -S_{233/002}+S_{223/001} -2S_{222}=
S_{024}+S_{123} \, .
\end{eqnarray*}

In short, the Specht basis and K-L basis for $n=3$  are

$$
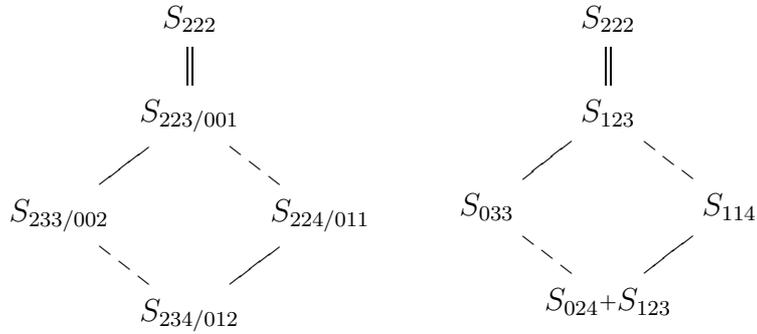
\xymatrix@R=0.5cm@C=0cm{
 &  &  & *{\GrTeXBox{S_{222} }}\arx2[d]& \\
 &  &  & *{\GrTeXBox{ S_{223/001}  }}\arx1[ld]\arx3[rd]& \\
 &  & *{\GrTeXBox{S_{233/002} }}\arx3[rd]&  
& *{\GrTeXBox{ S_{224/011} }}\arx1[ld]& \\
 &  &  & *{\GrTeXBox{ S_{234/012} }}& \\
} \ 
\xymatrix@R=0.5cm@C=0cm{
 &  &  & *{\GrTeXBox{ S_{222}  }}\arx2[d]& \\
 &  &  & *{\GrTeXBox{S_{123}}}\arx1[ld]\arx3[rd]& \\
 &  & *{\GrTeXBox{S_{033}}}\arx3[rd]&  
  & *{\GrTeXBox{ S_{114}}}\arx1[ld]& \\
 &  &  & *{\GrTeXBox{ S_{024}\plus S_{123} }}& \\
}
$$  

It seems a problem of interest for combinatorists to
give the explicit expression of $\KL_\l^S$ in terms of
Schur functions. The following lemma describes the case
where $\KL_\l^S$ coincides with a single Schur function.

\begin{lemma}
 Let $\l=[\beta^\alpha]$, $\beta\plus \alpha\leq n$
be a rectangular partition. Then
\begin{equation}
 \KL_\l^S = S_{ (n-1-\alpha)^\beta ,\, (n-1)^{n-\alpha-\beta},
  \, (n-1+\beta)^\alpha} \, .
\end{equation}
\end{lemma}
\proof  In the case of a rectangular partition, the
\emph{Kazhdan-Lusztig polynomials} 
are trivial (i.e. equal to $1$) \cite{Torun}. 
In our terms, this translates into the fact that 
$\KL_\l^S$ is the alternating sum of the elements of the Specht
basis on the interval of partitions contained in $\l$~:
$$ \KL_\l^S =  \sum_{\mu\le \l}
  (-1)^{|\mu|} S_{ ( (n-1)^n+\mu)/ \mu^\sim}  \, .$$
To compute this sum, one may suppose
that $h_i=h_i(\x)$, with $\x$ of cardinality $n$.
The skew Schur functions in the RHS factorize
into $S_{  (n-1)^n/ \mu^\sim}(\x) S_\mu(\x)$, 
according to (\ref{FactorSkewSchur}).

To avoid elaborate manipulations of determinants,
let us use the operator $\pi_\omega$ which sends
$x^v:\, v\in\N^n$ onto $S_{v_n,\dots, v_1}(\x)$.
One can now rewrite the RHS into
\begin{multline*}
\sum_{\nu\leq \alpha^\beta} (-1)^{|\nu|}
 x^{\overbrace{ \s n\moins 1,\dots,n\moins 1}^{n-\beta},
  n\moins 1\moins\nu_1,\dots, n\moins 1\moins\nu_\beta}
  \, S_{\nu^\sim}(\x) \, \pi_\omega  \\
  = x^{\overbrace{\s n\moins 1,\dots,n\moins 1}^{n-\beta},
   \overbrace{\s n\moins 1\moins\alpha,\dots,
    n\moins 1\moins\alpha}^{\beta}  } 
 S_{\beta^\alpha}(\x - x_{n-\beta+1}-\dots -x_n)\, \pi_\omega  \\
= x^{(n\moins 1)^{n-\beta}, (n\moins 1\moins\alpha)^\beta}
  x^{\beta^\alpha, 0^{n-\alpha}} \, \pi_\omega
= S_{(n-1-\alpha)^\beta ,\, (n-1)^{n-\alpha-\beta},
  \, (n-1+\beta)^\alpha}(\x) \, .
\end{multline*} 
This is the required identity.  \QED  

\subsection{Basis $\KL_\l^x$}

One generates it from $\KL_0^x=x_1\dots x_n$,
using the diagrams $\cD_\l$, 
interpreting an entry $[i,k]$
as $s_i + (1\plus k)^{-1}$.

$$ \xymatrix@R=0.5cm@C=-1.6cm{
 &  &  & *{\GrTeXBox{
              \KL_0^x=x^{111000} 
}}\arx2[d]& \\
 &  &  & *{\GrTeXBox{
          \KL_1^x= x^{110100}+x^{111000}
}}\arx1[ld]\arx3[rd]& \\
 &  & *{\GrTeXBox{
     \KL_{11}^x= x^{1011} \plus x^{1101}\plus x^{111}
}}\arx3[rd]&  & *{\GrTeXBox{
           \KL_2^x= x^{11001} \plus x^{1101}\plus x^{111}
}}\arx1[ld]& \\
 &  &  & *{\GrTeXBox{
  \begin{array}{c}
  \KL_{21}^x= x^{10101}+ x^{11001} \\
  +x^{1011} +x^{1101}+2 x^{111}  \end{array}
}}& \\
}$$

The coefficients are specializations $t=1$ of some
Kazhdan-Lusztig polynomials, which are, in the case
of Gra{\ss}mannians, easy to compute \cite{Torun}.
In the preceding example, there is only one non trivial 
Kazhdan-Lusztig polynomial, and it is 
equal to $1\plus t$.
This explains the coefficient $2$ in the expansion of  $\KL_{21}^x$.

\subsection{Dual basis  $\LK^a_\l$}

One generates it using the reversed graph,
with edges $s_i-1$, starting from
$$  \LK^a_\rho = \Delta^a(12|34| \dots |2n\moins 1,2n)
:= (a_1\moins a_2)(a_3\moins a_4) \dots (a_{2n-1} \moins a_{2n})\, .$$

Notice that $\LK^a_\rho$ is equal to the weight $\varphi^a(\rho)$.
This equality transfers in fact to all partitions.

\begin{lemma}  \label{Lemme LKa}
 For any $\l\leq \rho$, one has $ \LK^a_\l = \varphi^a(\l)$.
\end{lemma}
\proof The recursive definition of $ \LK^a_\l$
implies steps of the type
$$  \Delta^a( \dots | j,i |i\plus 1,k|\dots) 
 \hfl{s_i-1} 
 \Delta^a( \dots | j,i\plus 1 |i,k|\dots) -
 \Delta^a( \dots | j,i |i\plus 1,k|\dots)  \, .  $$ 
But, thanks to the Pl\"ucker relations for minors
of order $2$, this last element is equal to
$\Delta^a( \dots | j,k |i,i\plus 1|\dots)$. 
Therefore, the required property is true by 
decreasing induction on $\l$. \QED

\bigskip
$\xymatrix@R=0.5cm@C=0cm{
 &  &  & *{\GrTeXBox{
    \LK^a_{21}= \Delta^a(12|34|56)
}}\arx3[ld]\arx1[rd]& \\
 &  & *{\GrTeXBox{
     \LK^a_{11}= \Delta^a(12|36|45)
}}\arx1[rd]&  & *{\GrTeXBox{
    \LK^a_{2}= \Delta^a(14|23|56) 
}}\arx3[ld]& \\
 &  &  & *{\GrTeXBox{
  \LK^a_{1}= \Delta^a(16|23|45)
}}\arx2[d]& \\
 &  &  & *{\GrTeXBox{
  \LK^a_{0}= \Delta^a(16|25|34)
}}& \\
}$

\subsection{Dual basis  $\LK^x_\l$}

One defines elements in the ring 
$\cH_{2n}=\Z[x_1,\ldots, x_{2n}]/\Sym_+$,
where $\Sym_+$ is the ideal generated by symmetric
polynomials without constant term.
On this ring, one has a non-degenerate quadratic form
$(\, ,\, )^\d$ such that 
$(x^u\, ,\,x^v )^\d =(\moins 1)^{\ell(\sigma)}$ if there
exists a permutation $\sigma$ such that 
$(u+v)^\sigma = [2n\moins 1,\dots, 1,0]$,
and $(x^u\, ,\,x^v )^\d =0$  otherwise \cite{CBMS}.

One generates in $\cH_{2n}$  a family $\LK^x_\l$, $\l\leq \rho$,
using the reversed graph,
with edges $\, \moins(s_i\plus 1)$, starting from
$\LK^x_\rho= x^{0022\dots, 2n-2,2n-2}$.

Thanks to the ideal, the dual basis can be represented 
by single monomials~:
\begin{lemma} \label{Lemme LKx}
  For any $\l\leq \rho$, one has
   $ \LK^x_\l =  x^{\psi(\l)} $.
\end{lemma}
\proof
An elementary step $\LK^x_\l \to \LK^x_\mu$ 
in the recursive definition corresponds
to suppressing a corner labelled $[i,0]$ in $\cD_\l$,
and for the corresponding weight $\psi(\l)$, to the transformation
\begin{eqnarray*}
\psi(\l) &=&  w\, a\,  w'\,  \boxed{ab} \, w" b w'''   \\
\psi(\mu) &=&  w\, a\,  w' \, \boxed{bb} \,  w" a w'''   \\
\end{eqnarray*}
where the box stands in position $i,i\plus 1$ and $a,b$ are
two integers.
We claim that 
$$ x^{\dots \,\boxed{\s bb}\, \dots a\dots} 
\equiv - x^{\dots\,  \boxed{\s ab}\, \dots b\dots}(s_i\plus a) 
 = -\left( x^{\dots\, \boxed{\s ba}\, \dots b\dots} 
+ x^{\dots\, \boxed{\s ab}\, \dots b\dots}  \right)  \, ,$$
taking the notational liberty of replacing the common 
components of the two vectors by dots.
By permutation, one can shift the varying components
to the first three positions, and the property to show becomes
$$  x^{bba\dots} + x^{bab\dots} +x^{abb\dots}  \equiv 0 \, $$
the three monomials differing only in the exponents of
$x_1,x_2,x_3$.  
The nullity of the sum of three monomials can be tested
by checking the scalar products with all monomials $x^v$.
To hope for a permutation of $[2n\moins 1,\dots,0]$, the exponent
$v$ must belong to $\{0,1\}^{2n}$.
But in that case the polynomial 
$(x^{bba\dots} + x^{bab\dots} +x^{abb\dots}) x^v$
has at least a symmetry in $x_1,x_2$, or $x_1,x_3$,
or $x_2,x_3$ and therefore 
$(x^{bba\dots} \plus x^{bab\dots} \plus x^{abb\dots}, x^v)^\d$
is null even for those $v$.   \QED  

\bigskip
$$\xymatrix@R=0.5cm@C=0cm{
 &  &  & *{\GrTeXBox{
    \LK^x_{21}= x^{002244}
}}\arx3[ld]\arx1[rd]& \\
 &  & *{\GrTeXBox{
     \LK^x_{11}= x^{002442}
}}\arx1[rd]&  & *{\GrTeXBox{
    \LK^x_{2}= x^{022044}
}}\arx3[ld]& \\
 &  &  & *{\GrTeXBox{
  \LK^x_{1}= x^{022440}
}}\arx2[d]& \\
 &  &  & *{\GrTeXBox{
  \LK^x_{0}=x^{024420} 
}}& \\
}$$

\subsection{Duality $\KL_\l^\Delta$, $\LK_\l^x$}

It remains to justify the terminology \lq\lq dual basis'',
which does not reduce to reversing graphs.
It is natural to use the vanishing properties of Vandermonde determinants,
and, thus, to specialize the polynomials $\KL_\l^\Delta$.
In the present case, the following proposition
shows that $\psi(0),\ldots, \psi(\rho)$ are
convenient interpolation points.

\begin{proposition}
For any $\l,\mu\leq \rho$, one has
\begin{equation}
  \KL_\l^\Delta(\psi(\mu)) = (\moins 1)^{|\l|-|\mu|} c_n
                 \delta_{\l,\mu}  \, ,
\end{equation}
with $c_n = \prod_{1\leq i <j\leq n} (2i-2j)^2$.  
\end{proposition}

Let us denote $(f,x^v)$ the evaluation of a function 
$f(x_1,\dots,x_{2n})$ in $x_1=v_1,\dots, x_{2n}=v_{2n}$.
This form is compatible with the action of the symmetric
group: $(f s_i, g) =(f,g s_i)$. 

On the other hand, given three exponents which differ only
in three places, of the type $[\dots bb\dots a\dots],
 [\dots ba\dots b\dots],[\dots ab\dots b\dots]$,
one has
$$\bigl( \Delta^x(1\dots n|n\plus 1\dots 2n)\, ,\,
 x^{\dots bb\dots a\dots}+x^{\dots ba\dots b\dots}
+x^{\dots ab\dots b\dots}\bigr) =0 $$
and therefore, for any $\l\leq \rho$,
\begin{equation}  \label{Nul1}
 \bigl( \KL_\l^\Delta\, ,\, x^{\dots bb\dots a\dots}
   +x^{\dots ba\dots b\dots}
   +x^{\dots ab\dots b\dots}\bigr) =0   \, .
\end{equation}

Starting with 
$\bigl( \KL_\l^\Delta\, ,\, x^{002244\dots} \bigr) 
= c_n\delta_{\l,\rho}$,
one supposes that for some $\mu$ the proposition
is true. Let $[i,0]$ be a corner of $\cD_\mu$,
and let $\nu$ be the partition obtained from $\mu$
by removing this corner.  Then
$$  (\KL^x_\l, \LK^x_\nu)= 
   \bigl(\KL^x_\l, \moins \LK^x_\mu (s_i\plus 1) \bigr)
= \moins \bigl(\KL^x_\l (s_i\plus 1), \LK^x_\mu) \, $$
thanks to (\ref{Nul1}) and Lemma \ref{Lemme LKx}.
Non nullity can occur only for $\l=\mu$ or $\l=\nu$.
Since $\KL^x_\nu (s_i\moins1) = \KL^x_\mu 
  +\sum_{\eta\neq \mu,\nu} c_\eta \KL^x_\eta$, one has
\begin{eqnarray*}
 \moins \bigl(\KL^x_\nu (s_i\plus 1), \LK^x_\mu) 
 &=& \moins \bigl(\KL^x_\mu, \LK^x_\mu)  \\ 
 \moins \bigl(\KL^x_\mu (s_i\plus 1), \LK^x_\mu) 
 &=&  \moins \bigl(\KL^x_\nu (s_i\moins1)(s_i\plus 1),\LK^x_\mu) 
     =0  \, ,
\end{eqnarray*}
which is what is expected for the proof by induction to be valid.  
\QED

\subsection{Duality $\KL_\l^x$, $\LK_\l^x$}
The polynomials $\KL_\l^x \LK_\mu^x$ have total degree 
$0\plus 1\plus \cdots\plus (2n \moins 1)$. This points to
using the form $(\, ,\, )^\d$ \cite{Turbo}.

\begin{proposition}
For any $\l,\mu\leq \rho$, one has
\begin{equation}   \label{DualiteKL}
  \bigl( \KL_\l^x, \LK_\mu\bigr)^\d= \delta_{\l,\mu}  \, .
\end{equation}
\end{proposition}
\proof  The starting point is 
$ \bigl( x^{1\dots 1 0\dots 0}, \LK_\mu\bigr)^\d= \delta_{0,\mu}$.
The general case is deduced by the same induction as 
in the preceding case, using, for any $f,g$, any $i\leq 2n\moins 1$
the identity $(f,gs_i)^\d = - (f s_i,g)^\d$.   \QED

\subsection{Duality $\KL_\l^x$, $\LK_\l^a$}

This time, we shall use the vanishing properties of 
the polynomials $\LK_\l^a$.
\begin{proposition}
For any $\l,\mu\leq \rho$, one has
\begin{equation}
  ( \LK_\l^a, \KL_\mu^x) = \delta_{\l,\mu} \, .
\end{equation}
\end{proposition}
\proof  The starting point is
$\LK_\mu^a\left(x^{1\dots 1 0\dots 0} \right)= \delta_{0,\mu}$
and the proof by induction goes as before.  \QED

\section{Hankel Pfaffians in terms of Kazhdan-Lusztig basis}

In \cite[Th.4.1]{LaPfaff}, we have given several expressions
of a Pfaffian with entries $(a_i\moins a_j) g_{i,j}$.

For the present case, for $n=3$, this would read in the
Specht basis as
\begin{multline*}
\Pf(\a,h(\x),3,0)  = 
 -\smallyoung{4 &5&6\cr 1&2&3\cr}\, S_{444}(\x)
 +\smallyoung{3&5&6\cr 1&2&4\cr}\, S_{445/001}(\x)
 - \smallyoung{3&4&6\cr 1&2&5\cr}\, S_{446/011}(\x)\\
 -\smallyoung{2&5&6\cr 1&3&4\cr}\, S_{455/002}(\x)
 +\smallyoung{2&4&6\cr 1&3&5\cr}\, \Bigl(S_{456/012}(\x)
  -S_{444}(\x)\Bigr) \, ,
\end{multline*}
each tableau being interpreted as as a product of
factors $(a_i\moins a_j)$ corresponding to its columns.

In \cite{LaPfaff}, we have shown in particular that the Pfaffian is diagonal 
in terms of Young's orthonormal basis. The underlying
quadratic form in that case is formally defined
in terms of tableaux \cite{Rutherford,Turbo},
but corresponds to the form $(\, ,\, )^\d$
when interpreted in the appropriate spaces.
Therefore, thanks to (\ref{DualiteKL}),
the Pfaffian remains diagonal when using 
Kazhdan-Lusztig bases instead of Young's bases \cite{Young},
and one has the following theorem.

\begin{theorem}     \label{th:PfaffKL}
Given $n$, then one has
\begin{equation}   \label{Pfaff2KL}
 \Pf(\a,\h,n,1\moins n)  =
  \sum_{\l\le \rho}  (\moins 1)^{|\l|} \LK^a_\l \, \KL^S_\l  \, .
\end{equation}
\end{theorem}

For example, for $n=3$, one has
\begin{multline*}
 \Pf(\a,\h,3,-2) =
 \Delta^a(16|25|34) \KL^S_0
 -\Delta^a(16|23|45) \KL^S_1
  +\Delta^a(14|23|56) \KL^S_2   \\
  +\Delta^a(12|36|45) \KL^S_{11}
  -\Delta^a(12|34|56)\KL^S_{21}  \\
 =                    
   (a_1 - a_6) (a_2 - a_5) (a_3 - a_4)\, S_{222}
- (a_1 - a_6) (a_2 - a_3) (a_4 - a_5)\, S_{123} \\  
+ (a_1 - a_4) (a_2 - a_3) (a_5 - a_6)\, S_{114}  
+ (a_1 - a_2) (a_3 - a_6) (a_4 - a_5)\, S_{033}  \\
- (a_1 - a_2) (a_3 - a_4) (a_5 - a_6)(S_{024}+S_{123}) \, .
\end{multline*}

\section{Discriminants and Bezoutians}

The case where $\a$ specializes to $\q=[1,q,\dots,q^{2n-1}]$
is of special interest. 
The Pfaffian $\Pf(\q, \h(\x\moins \y,k)$ is proportional 
$\Pf(\q, \h(\x),2-n)$, which is equal to the determinant
of the matrix $\M_h(\q,\x,0)$.
But, because of homogeneity, one does not change the value of 
this determinant
by replacing the entries $ a_i h_j(\x) = q^{i-1} h_j(\x)$ by 
$q^j h_j(\x) =  h_j(q\x)$. The Laplace expansion of this new matrix
along the first $n$ columns is equal to $q^{\binom{n}{2}}$
times the expansion of $S_{n^n}(\x -q\x)= R(\x,q\x)=
 (1\moins q)^n (x_1\dots x_n) \cD(\x,q)$.

In final, one has the following theorem linking Pfaffians,
resultants and discriminants.

\begin{theorem}    \label{th:PfaffDiscri}
Given $\x$ of cardinality $n$, $\y$ of cardinality $m$,
and $k\in\Z$, then 
\begin{equation}    \label{PfaffDiscri}
\Pf(\q, \h(\x\moins \y),k)  =     q^{\binom{n}{2}} 
 (1\moins q)^n (x_1\dots x_n)^{k+n-1-m} R(\x,\y)\, \cD(\x,q)   \, .
\end{equation}
\end{theorem}

Consequently, one can use the expression of the Pfaffian in terms
of the KL-basis to expand the $q$-discriminant.

For example, for $n=2,3,4$, denoting $[i]$ the $q$-integer
$(q^i\moins 1)/(q\moins 1)$, one has the following
expansions~:

$$\cD(2,q) = [3]  \KL^S_0 - q\KL^S_1 =(1\plus q\plus q^2)S_{11}
    - q S_{02} \, , $$
\begin{multline*}
\cD(3,q) = [3] [5] \KL^S_0 -q [5] \KL^S_1 +q [3]^2\KL^S_2
                    + q^2[3] \KL^S_{11}  - q^3  \KL^S_{12}   \\
= (1\plus q\plus q^2)(1\plus\dots\plus q^4) S_{222}
  -q(1\plus\dots\plus q^4)S_{123} +
  q^2(1\plus q\plus q^2)S_{114}   \\
 +  q^2(1\plus q\plus q^2)S_{033} - q^3 (S_{024}+S_{123})\, ,
\end{multline*}
\begin{multline*}
\cD(4,q)=
 [3]  [5]  [7] \KL^S_0 -q  [5]  [7] \KL^S_1 +q^2  [3]  [7] \KL^S_2 
        -q^3  [3]  [5] \KL^S_3\\
+q^2  [3]  [7] \KL^S_{11}  -q^3  [7] \KL^S_{21} -q^3  [3]  [5] \KL^S_{111}   \\
 +q^4  [3]^2 \KL^S_{22} 
 + q^4  [5] \KL^S_{31} +q^4  [5] \KL^S_{211} \\  
 -q^5  [3] \KL^S_{32}
-q^5  [3] \KL^S_{221}
   -q^5  [3] \KL^S_{311} +q^6 \KL^S_{321} \, .
\end{multline*}

Since one has $p_r(\x) = n h_r(\x-\x^{der})$, 
one obtains from (\ref{PfaffDiscri}) the evaluation of
Pfaffians where the power sums replace the complete functions,
the resultant being replaced by the discriminant.

\begin{corollary}
Given $\x$ of cardinality $n$ and $k\in\Z$, then

 \begin{equation}    \label{PfaffDiscri_p}
\Pf(\q\, , {\mathbf p}(\x),k)  =
  (\moins q)^{\binom{n}{2}}    
   (1\moins q)^n (x_1\dots x_n)^{k}\,  \cD(\x,1)\, \cD(\x,q)   \, .
\end{equation}

\end{corollary}

The limit $q\to 1$ gives that the Pfaffian with entries
$(i-j) p_{i+j-3+k}(\x)$ is equal to
$$  (x_1\dots x_n)^{k+n-1-m}\,  \cD(\x,1)\cD(\x,q) $$
and that the Pfaffian with entries
$(i-j) p_{i+j-3+k}(\x)$ is equal to
$$  (\moins 1)^{\binom{n}{2}}
    (x_1\dots x_n)^{k}\,  \cD(\x,1)^2   \, .$$

One can also use the matrix $\M_h(\q,\x -\x^{der},k)$.
For example, for $n=2$, $k=0$, one has
$$
\begin{vmatrix}
x_{1}^{2}\plus x_{2}^{2}   &
x_{1}^{3}\plus x_{2}^{3}&
q^{3}(x_{1}^{2}\plus x_{2}^{2})&
q^{3}(x_{1}^{3}\plus x_{2}^{3}) \\
x_{1}\plus x_{2} &
x_{1}^{2}\plus x_{2}^{2}&
q^{2}(x_{1}\plus x_{2})&
q^{2}(x_{1}^{2}\plus x_{2}^{2}) \\
2 &
x_{1}\plus x_{2}&
2q &
q(x_{1}\plus x_{2})\\
\frac{1}{x_1}\plus \frac{1}{x_2}&
2&
\frac{1}{x_1}\plus \frac{1}{x_2}&
2
\end{vmatrix}    = \frac{q (q-1)^2}{x_1x_2}
(x_1 \moins x_2)^4 (x_1\moins qx_2)(x_2\moins qx_1) \, .
$$

\section{Remark about Macdonald polynomials}

Tke KL-basis for the representations of shape $[2^n]$ or $[n,n]$
of the Hecke algebra is related to the \emph{non-symmetric
Macdonald polynomials}  $M_v(\x_{2n};t,q)$, $v\in\N^{2n}$.
In \cite{deGierLa}, one finds a common deformation of
the KL basis and the Macdonald polynomials indexed
by a permutation of $[\dots 22 11 00]$, when $q$ is specialized
to a certain root of $t$.

The $t$-discriminant itself, which is a symmetric function,
is equal to the specialization of the symmetric 
Macdonald polynomial indexed by the (decreasing) partition
$[2n\moins 2,\dots, 2,0]$ \cite[Th.3.2]{BouLuque}. The
$t$-discriminants also appear as specializations of symmetric
Macdonald polynomial indexed by rectangular partitions
\cite[Remark 4.9]{Luque}. 

Investigating extensively the specializations of symmetric or 
non-symmetric Macdonald polynomials
at $q=t^\alpha$ seems to be of great interest.


\end{document}